\documentclass[preprint]{elsarticle}

\usepackage{amsmath}
\usepackage{amssymb}
\usepackage{graphicx}

\DeclareMathOperator{\Gr}{Gr}
\newcommand{\flip}[1]{#1^\Diamond}

\newcommand{\noflip}[1]{#1^{\phantom{\Diamond}}}
\newcommand{\x}{{\scriptstyle \times}}
\newcommand{\z}{{\scriptstyle 0}}
\newcommand{\one}{{\scriptstyle 1}}

\newcommand{\degree}{\ensuremath{^\circ}}

\newtheorem{theorem}{Theorem}
\newproof{proof}{Proof}

\journal{Linear Algebra and its Applications}
\begin{document}
\begin{frontmatter}

\title{Banded Householder representation of linear subspaces}
\author{G. Irving}
\ead{irving@naml.us}
\address{Weta Digital, 9-11 Manuka St., Miramar, Wellington, New Zealand}

\begin{abstract}
We show how to compactly represent any $n$-dimensional subspace of $\mathbf{R}^m$ as a banded product
of Householder reflections using $n(m-n)$ floating point numbers.  This is optimal since these subspaces form a
Grassmannian space $\Gr_n(m)$ of dimension $n(m-n)$.
The representation is stable and easy to compute: any matrix can be factored into the product
of a banded Householder matrix and a square matrix using two to three QR decompositions.
\end{abstract}



\end{frontmatter}

If $m \ge n$, the Householder QR algorithm represents an $m \times n$ orthogonal matrix $U$ as a product of $n$ Householder reflections using
a total of $n(m-(n+1)/2)$ floating point numbers \cite[chap.\ 5]{golub1996matrix}.  However, in some applications only the range of $U$ is important; any other orthogonal
matrix with the same range is equivalent.  One example is the hierarchically semiseparable representations of \cite{hss2004}, where a tree
of orthogonal matrices is used to compress matrices with significant offdiagonal structure.  Since the $n$-dimensional subspaces
of $\mathbf{R}^m$ form a Grassmannian manifold $\Gr_n(m)$ of dimension $n(m-n)$, we expect that \emph{some} orthogonal matrix with the correct
range can be represented with $n(m-n)$ floats.  The following theorem provides such a representation:

\begin{theorem}\label{banded}
If $m \ge n$, any $A \in \mathbf{R}^{m \times n}$ can be factored as
\begin{equation*}
A = G \left(\begin{array}{c}B \\ 0\end{array}\right)
\end{equation*}
where $B \in \mathbf{R}^{n \times n}$ is
square and $G$ is a product of $n$ Householder reflections with banded structure:
\begin{equation}\label{bh}
\begin{array}{c}
\begin{aligned}
G &= H_1 H_2 \cdots H_n \\
H_i &= I - \frac{2 v_i v_i^T}{v_i^T v_i}
\end{aligned} \\
v_i = \left(\underbrace{0, 0, \cdots, 0}_{i-1}, 1, \underbrace{v_{i,i+1}, \cdots, v_{i,i+m-n}}_{m-n}, \underbrace{0, 0, \cdots, 0}_{n-i} \right)^T \\
\end{array}
\end{equation}
Since each $v_i$ has $m-n+1$ nonzero components, the first of which is $1$, $G$ can be stored in $n(m-n)$ floats.  If $A$ is full rank,
the matrices $G$ and $B$ are unique, although the Householder vectors $v_i$ may not be.
\end{theorem}
\begin{proof}
Observe that (\ref{bh}) is exactly the factored form produced by standard Householder QR except for the trailing zeroes in each $v_i$, which correspond to the
extreme lower triangle $i>j+m-n$.  To introduce these zeroes, define $\flip{A}$ as $A$ rotated by 180\degree~(pronounced ``flip A''),
\begin{equation*}
\flip{A}_{ij} = A_{m-i+1,n-i+1}
\end{equation*}
and perform an LQ decomposition of $\flip{A}$:
\begin{align*}
\flip{A} &= \noflip{L} \noflip{Q} \\
\noflip{A} &= \flip{L} \flip{Q}
\end{align*}
Since $L_{ij} = 0$ for $i<j$, $\flip{L}_{ij} = 0$ for $i>j+m-n$.  The Householder QR algorithm constructs $v_i$ as a linear combination of $e_i$
and the $i$th column of the matrix (after rotation by the previous Householder reflections), and the first component of each vector can be chosen
to be 1 \cite[chap. 5]{golub1996matrix}.
Therefore, a Householder QR decomposition
\begin{equation*}
\flip{L} = G \left(\begin{array}{c}R \\ 0\end{array}\right)
\end{equation*}
will produce $G$ with the correct banded structure.  Our final factorization is
\begin{equation*}
A = G \left(\begin{array}{c}R \\ 0\end{array}\right) \flip{Q} = G \left(\begin{array}{c}B \\ 0\end{array}\right)
\end{equation*}
The steps are visualized in Figure~\ref{sq}.  Uniqueness of $G$ and $B$ follows from the uniqueness of the QR decomposition when $A$ is full
rank \cite[chap. 5]{golub1996matrix}.  Note that $B$ is orthogonal whenever $A$ is orthogonal.
\qed\end{proof}

\begin{figure}
\begin{center}
{
\renewcommand{\tabcolsep}{0cm}
\renewcommand{\arraystretch}{1}
\begin{equation*}
\begin{array}{c@{}c@{}c@{}c@{}c@{}c@{}c}
{\scriptscriptstyle
\left(\begin{array}{cccc}
\x & \x & \x & \x \\
\x & \x & \x & \x \\
\x & \x & \x & \x \\
\x & \x & \x & \x \\
\x & \x & \x & \x \\
\x & \x & \x & \x \\
\x & \x & \x & \x \\
\x & \x & \x & \x
\end{array}\right)}
& \Rightarrow &
{\scriptscriptstyle
\left(\begin{array}{cccc}
\x & \z & \z & \z \\
\x & \x & \z & \z \\
\x & \x & \x & \z \\
\x & \x & \x & \x \\
\x & \x & \x & \x \\
\x & \x & \x & \x \\
\x & \x & \x & \x \\
\x & \x & \x & \x
\end{array}\right)}
& \Rightarrow &
{\scriptscriptstyle
\left(\begin{array}{cccc}
\x & \x & \x & \x \\
\x & \x & \x & \x \\
\x & \x & \x & \x \\
\x & \x & \x & \x \\
\x & \x & \x & \x \\
\z & \x & \x & \x \\
\z & \z & \x & \x \\
\z & \z & \z & \x
\end{array}\right)}
& \Rightarrow &
{\scriptscriptstyle
\left(\begin{array}{cccc}
\one & \z & \z & \z \\
\x & \one & \z & \z \\
\x & \x & \one & \z \\
\x & \x & \x & \one \\
\x & \x & \x & \x \\
\z & \x & \x & \x \\
\z & \z & \x & \x \\
\z & \z & \z & \x
\end{array}\right)} \\
\multicolumn{7}{c}{\vspace{-.2cm}} \\
A & &
L = \flip{A} Q^T & &
\flip{L} & &
\left(\begin{array}{ccc} v_1 & \cdots & v_n \end{array} \right)
\end{array}
\end{equation*}
}
\end{center} \caption{\label{sq} To compute a banded Householder factorization of a matrix $A$, we perform an LQ factorization of $A$ rotated by
180\degree~(denoted $\flip{A}$) to zero the extreme lower triangle of $A$, then perform a QR decomposition to zero the upper triangle and
construct the banded Householder vectors $v_i$.}
\end{figure}

Since the construction uses only matrix multiply and Householder QR decomposition as primitives, the computation of $G$ is stable and requires $O(mn^2)$ flops.
Normally the scalar factors $\beta_i = 2/(v_i^T v_i)$ will be precomputed and stored for an additional $n$ floats of storage.  Once $\beta_i$
is available, a single Householder reflection requires $4(m-n)+2$ flops and matrix vector products $G x$ or $G^T y$ can be computed in $4n(m-n)+2n$
flops.  If blocking is desired, we can represent products of $b$ Householder
reflections where $b$ is the block size using the $I - V T V^T$ representation \cite{schreiber1987wy}, which involves a relatively small increase in
storage and flops if $b \ll m-n$.

Unfortunately, the banded Householder representation in Theorem~\ref{banded} is inefficient for large subspaces of $\mathbf{R}^m$, since $G x$
involves a large number of small level 1 BLAS operations if $m-n$ is small.  To remedy this problem, we can represent a large $n$-dimensional
subspace in terms of its small $(m-n)$-dimensional orthogonal complement as follows:

\begin{theorem}\label{upsidedown}
If $m \ge n$, any $A \in \mathbf{R}^{m \times n}$ can be factored as 
\begin{equation*}
A = G \left(\begin{array}{c}0 \\ B\end{array}\right)
\end{equation*}
where $B \in \mathbf{R}^{n \times n}$
is square and $G$ is a banded product of $m-n$ Householder reflections with $n+1$ nonzero components per vector,
the first of which is 1 (equivalent to (\ref{bh}) with $m-n$ and $n$ swapped).  In particular, $G$ can also be
stored in $n(m-n)$ floats.
\end{theorem}
\begin{proof}
Perform a QR decomposition of $A$ to get
\begin{align*}
A &= \left(\begin{array}{cc}U_1 & U_2\end{array}\right)\left(\begin{array}{c}R \\ 0\end{array}\right) = U_1 R
\end{align*}
Here the column span of $U_1 \in \mathbf{R}^{m \times n}$ contains the range of $A$, and the span of $U_2 \in \mathbf{R}^{m \times (m-n)}$ is contained
inside the nullspace of $A^T$.  Banded Householder factorization of $U_2$ gives
\begin{align*}
U_2 &= G \left(\begin{array}{c}Q \\ 0 \end{array}\right)
     = \left(\begin{array}{cc} G_1 & G_2\end{array}\right) \left(\begin{array}{c}Q \\ 0 \end{array}\right)
     = G_1 Q
\end{align*}
whence $G_1 = U_2 Q^T$ and
\begin{align*}
G^T A &= \left(\begin{array}{c} G_1^T \\ G_2^T \end{array}\right) A
       = \left(\begin{array}{c} Q U_2^T A \\ G_2^T A \end{array}\right)
       = \left(\begin{array}{c} 0 \\ G_2^T A \end{array}\right)
\end{align*}
since $U_2^T A = U_2^T U_1 R = 0$.  Our final decomposition is
\begin{equation*}
A = G \left(\begin{array}{c}0 \\ G_2^T A \end{array}\right) = G \left(\begin{array}{c}0 \\ B\end{array}\right)
\end{equation*}
\qed\end{proof}

Using Theorem~\ref{banded} for $m-n \ge n$ and Theorem~\ref{upsidedown} for $m-n < n$, the resulting $G$ consists of at most $m/2$ Householder vectors
each with at least $m/2+1$ nonzero components.  In particular, the blocked $I - V T V^T$ form is efficient whenever $b \ll m,n$, regardless
of the value of $m-n$.

\section{Application}

\begin{figure}
\begin{center}
\includegraphics[width=1.5in]{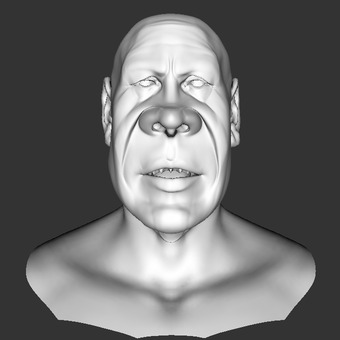}
\includegraphics[width=1.5in]{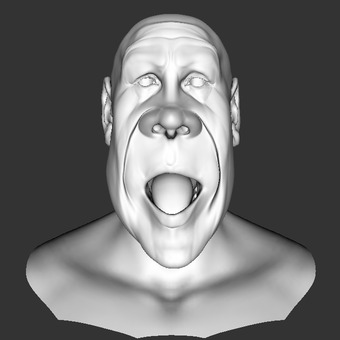}
\includegraphics[width=1.5in]{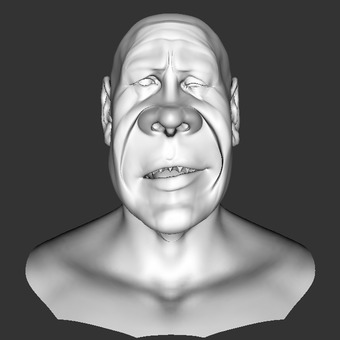}
\end{center}
\caption{\label{genman} Frames from an animation of a digital character using blend shapes stored in a hierarchically semiseparable
representation (HSS).  Using banded Householder form for the rotations in the HSS tree reduces the storage costs by 45.7\% over dense storage,
or 29.5\% of Householder storage.}
\end{figure}

Our motivating application for the banded Householder decomposition is the compression of blend shape matrices for digital characters.  We start with
a large, mostly dense matrix where each column represents a pose of the digital character mesh.  An example face with 42391 vertices and 730 blend
shapes is shown in Figure~\ref{genman}.  The original matrix consumes 348 MB of storage in single precision.  To reduce
this, we compute a lossy hierarchically semiseparable (HSS) representation for the matrix \cite{hss2004}, which represents a matrix as a tree of
rotations and dense blocks.  If the rotations are stored in dense form, the HSS representation requires 46.8 MB of storage.  Using Householder
form reduces the storage cost to 36.0 MB (77.7\% of dense), and banded Householder form reduces the cost further to 25.4 MB (54.3\% of dense).
On an 8 core Intel Xeon 2.8 GHz machine, the cost to multiply the HSS representation with a vector is 11.2 ms using dense storage with optimized
BLAS and 10.7 ms using banded Householder storage with handwritten, unvectorized C.  Since the required memory traffic in the banded Householder
case is roughly half that of the dense case, we expect this comparison would improve significantly if the banded Householder code were appropriately
vectorized.

\section{Acknowledgements}
I am grateful to Tamar Shinar, J.P.\ Lewis, and Jaewoo Seo for helpful discussions and comments on the paper.

\bibliography{references}
\bibliographystyle{model1-num-names}
\end{document}